\documentclass[12pt, reqno]{amsart}
\usepackage{amsfonts}
\usepackage{amssymb}
\usepackage{amsmath}

\newtheorem{theorem}{Theorem}
\theoremstyle{plain}

\newtheorem{corollary}[theorem]{Corollary}

\begin{document}
\title[] {Some identities on %On the $p$-adic integral representation
                             %on $\mathbb Z_p$ of the 
Bernstein polynomials associated with $q$-Euler polynomials}
\author{A. Bayad}
\address{Abdelmejid Bayad. D\'epartement de math\'ematiques \\
Universit\'e d'Evry Val d'Essonne, Bd. F. Mitterrand, 91025 Evry
Cedex, France  \\} \email{abayad@maths.univ-evry.fr}
\author{T. Kim}
\address{Taekyun Kim. Division of General Education-Mathematics \\
Kwangwoon University, Seoul 139-701, Republic of Korea  \\}
\email{tkkim@kw.ac.kr}
\author{B. Lee}
\address{Department of Wireless Communications Engineering,  \\
Kwangwoon University, Seoul 139-701, Republic of Korea  \\}
\email{bjlee@kw.ac.kr}

\author{S.-H. Rim}
\address{Department of Mathematics Education, \\
Kyungpook National University, Taegu 702-701, Republic of Korea  \\}
\email{shrim@knu.ac.kr}

%\thanks{$*$ correponding author}
\maketitle

{\footnotesize {\bf Abstract} \hspace{1mm} {
In this paper we investigate some interesting properties of the $q$-Euler polynomials. The purpose of this paper is to give some relationships between Bernstein and $q$-Euler polynomials which are derived by the $p$-adic integral representation of the Bernstein polynomials associated with $q$-Euler polynomials.
}

\section{Introduction}
Let $p$ be a fixed odd prime number. Throughout this paper %The symbol
$\Bbb Z_p$, $\Bbb Q_p$ and $\Bbb C_p$ denote the ring of $p$-adic
integers, the field of $p$-adic numbers and the field of $p$-adic
completion of the algebraic closure of $\Bbb Q_p$, respectively (see
[1--15]). Let $\Bbb N$ be the set of natural numbers and $\Bbb Z_+ =\Bbb N
\cup \{0\}$. 
The normalized $p$-adic absolute value is defined by $|p|_p=\frac1{p}.$ As an indeterminate, we assume that  $q\in \Bbb C_p$ with $|1-q|_p <1$.
Let $UD(\Bbb Z_p)$ be the space of  uniformly differentiable function  on $\Bbb Z_p$. For $f\in UD(\Bbb Z_p)$, the $p$-adic invariant integral on $\Bbb Z_p$ is defined by  
\begin{eqnarray}
I_{-1} (f)=\int_{\Bbb Z_p } f(x) d\mu_{-1} (x)&=&
\lim_{N \rightarrow \infty}  \sum_{x=0}^{p^N-1} f(x)\mu_{-1}(x+p^N\Bbb Z_p)\\
%\end{eqnarray}$$\qquad\qquad \qquad\qquad\qquad \qquad 
&=&
\lim_{N \rightarrow \infty}  \sum_{x=0}^{p^N-1} f(x)(-1)^x, \quad (\text{see \cite{5,6,7,8}}).\nonumber
\end{eqnarray}
For $n \in \Bbb N$, we can derive the following integral equation from (1): 
\begin{eqnarray}
I_{-1} (f_n)=(-1)^n \int_{\Bbb Z_p } f(x) d\mu_{-1} (x)%I_{-1}(f)
 + 2\sum_{l=0}^{n-1}(-1)^{n-1-l}f(l),
\end{eqnarray}
where $f_n (x)=f(x+n)$ ( see \cite{5,6,7,8,9}).
As well known definition, the Euler polynomials are given by the generating function as follows:
\begin{eqnarray}
\frac{2}{e^t +1}e^{xt} =e^{E(x)t}= \sum_{n=0}^{\infty} E_n (x) \frac{t}{n!},
\quad \text{(see \cite{5,6,7,8,9,10,11,12,13,14,15})},
\end{eqnarray}
with usual convention about replacing $E^n(x)$ by $E_n(x)$. 
In the special case $x=0$, $E_n(0)=E_n$ are called the $n$-th Euler
numbers. From (3), we can derive the following recurrence formula for Euler numbers
\begin{eqnarray*}
E_0=1, (E+1)^n+E_n=\left\{\begin{array}{ll} 2 \ \ &\hbox{if}\ \ n=0,
\vspace{2mm}\\
0 \ \
&\hbox{if}\ \ n>0,\quad(\text{see \cite{10}}).
\end{array}\right.  \notag
\end{eqnarray*}
with usual convention about replacing $E^n$ by $E_n$.
By the definitions of Euler numbers and polynomials, we get
\begin{eqnarray*}
E_n(x)=(E+x)^n=\sum_{l=0}^{n} \binom{n}{l}x^{n-l}E_l, \quad(\text{see \cite{5,6,7,8,9,10,11,12,13,14,15}}).
\end{eqnarray*}
Let $C[0, 1]$ denote the set of continuous functions on $[0, 1]$.
For $f \in C[0, 1]$, Bernstein introduced the following well-known
linear positive operator in the field of real numbers $\Bbb R$ :
\begin{eqnarray*}
\mathbb{B}_n(f | x)= \sum_{k=0}^{n} f(\frac{k}{n})\binom{n}{k} x^k
(1-x)^{n-k}=\sum_{k=0}^{n} f(\frac{k}{n})B_{k,n}(x),
\end{eqnarray*}
where
$\binom{n}{k}=\frac{n(n-1)\cdots(n-k+1)}{k!}=\frac{n!}{k!(n-k)!}$
(see \cite{1, 2, 5, 9, 10, 14}). Here, $\mathbb{B}_n (f|x)$ is called the
Bernstein operator of order $n$ for $f$. For $k,n \in \Bbb Z_+$, the Bernstein polynomials of degree $n$ are defined by
\begin{eqnarray}
B_{k,n}(x)= \binom{n}{k}x^k (1-x)^{n-k}, \quad \text{for} \ x\in [0,1].
\end{eqnarray}
In this paper, we study the properties of $q$-Euler numbers and
polynomials. From these properties we investigate some identities on the
$q$-Euler numbers and polynomials. Finally, we give some relationships between Bernstein and $q$-Euler polynomials, which are derived by the $p$-adic integral representation of the Bernstein polynomials associated with $q$-Euler polynomials.
\section{ $q$-Euler numbers and polynomials}
In this ection we assume that $q\in\Bbb C_p$ with $|1-q|_p<1.$  Let $f(x)=q^xe^{xt}.$ From (1) and (2), we have 
\begin{eqnarray}
\int_{\Bbb Z_p } f(x) d\mu_{-1} (x)=\frac{2}{qe^t +1}e^{xt}. 
\end{eqnarray}
Now, we define the $q$-Euler numbers as follows:
\begin{eqnarray}
\frac{2}{qe^t +1} =e^{E_qt}=\sum_{n=0}^{\infty} E_{n,q}\frac{t^n}{n!},
\end{eqnarray}
with the usual convention about replacing $E_q^n$ by $E_{n,q}.$ 

By (6), we easily get

\begin{eqnarray}
E_{0,q}=\frac{2}{q+1}, \text{ and }  q(E_q+1)^n+E_{n,q}=\left\{\begin{array}{ll} 2 \ \ &\hbox{if}\ \ n=0,
\vspace{2mm}\\
0 \ \
&\hbox{if}\ \ n>0,
\end{array}\right.  %\notag
\end{eqnarray}
with usual convention about replacing  $E_q^n$ by $E_{n,q}.$ \\
 We note that
 \begin{eqnarray}
\frac{2}{qe^t +1} =\frac{2}{e^t+q^{-1}}.\frac{2}{1+q}=\frac{2}{1+q}\sum_{n=0}^{\infty} H_n(-q^{-1})\frac{t^n}{n!},
\end{eqnarray}
where $H_n(-q^{-1})$ is the $n$-th Frobenius-Euler numbers.\\
From (5), (6) and (8), we have
\begin{eqnarray}
\int_{\Bbb Z_p } q^xe^{xt} d\mu_{-1} (x)=E_{n,q}=\frac{2}{1+q}H_n(-q^{-1}), \textrm{ for } n\in\Bbb Z_{+}.
\end{eqnarray}
Now, we consider the $q$-Euler polynomials as follows:
\begin{eqnarray} \frac{2}{qe^t +1}e^{xt} =e^{E_q(x)t}=\sum_{n=0}^{\infty} E_{n,q}(x)\frac{t^n}{n!},
\end{eqnarray}
whith the usual convention  $E_q^n(x)$ by $E_{n,q}(x).$ 

From (2), (5) and (10), we get 
\begin{eqnarray} 
\int_{\Bbb Z_p } q^xe^{(x+y)t} d\mu_{-1} (y)=\frac{2}{qe^t +1}e^{xt}=\sum_{n=0}^{\infty} E_{n,q}(x)\frac{t^n}{n!}.
\end{eqnarray}
By comparing the coefficients on the both sides of (10) and (11), we
get the following Witt's formula for the $q$-Euler polynomials as follows:
\begin{eqnarray}
\int_{\Bbb Z_p} q^y(x+y)^n d\mu_{-1}(y)=E_{n,q}(x)=\sum_{l=0}^{n}\binom{n}{l}x^{n-l}E_{l,q}.
\end{eqnarray}
From (10) and (12), we can derive the following equation:
\begin{eqnarray} 
\frac{2q}{qe^t +1}e^{(1-x)t} =\frac{2}{1+q^{-1}e^{-t}}e^{-xt}=\sum_{n=0}^{\infty} E_{n,q^{-1}}(x)(-1)^n\frac{t}{n!} .
\end{eqnarray}
By (10) and (13), we obtain the following reflection symmetric property for the $q$-Euler polynomials. 
\begin{theorem}For $n\in\Bbb Z_{+}$, we have 
$$(-1)^n E_{n,q^{-1}}(x) = qE_{n,q}(1-x).$$
\end{theorem}
From (9), (10), (11) and (12), we can derive the following equation: for $n\in\Bbb N$,
\begin{eqnarray*}
E_{n,q}(2)&=&(E_q+1+1)^n=\sum_{l=0}^{n}\binom{n}{l}E_{l,q}(1)\\
&=& E_{0,q}+\frac1{q}\sum_{l=1}^{n}\binom{n}{l}qE_{l,q}(1)= \frac2{1+q}- \frac1{q}\sum_{l=1}^{n}\binom{n}{l}E_{l,q}\\
&=& \frac2{1+q}+\frac2{q(1+q)}- \frac1{q}\sum_{l=0}^{n}\binom{n}{l}E_{l,q} \\
&=& \frac2{q}-
\frac1{q^2}qE_{n,q}(1)=\frac2{q}+\frac1{q^2}E_{n,q},~\textrm{ by using
  recurrence formula (7)}.
\end{eqnarray*}
%%%%%%%%%%%%%%%%%%%
Therefore, we obtain the following theorem.
\begin{theorem}For $n\in\Bbb N$, we have 
$$qE_{n,q}(2)=2+\frac1{q}E_{n,q}.$$
\end{theorem}
By using (9) and (12), we get 
\begin{eqnarray*}
\int_{\Bbb Z_p} q^{-x}(1-x)^n d\mu_{-1}(x)
&=& (-1)^n \int_{\Bbb Z_p}  q^{-x}(x-1)^n d\mu_{-1}(x) \notag \\
&=&(-1)^nE_{n,q^{-1}}(-1)=q\int_{\Bbb Z_p} (x+2)^n d\mu_{-1}(x)=q\left(\frac2{q}+\frac1{q^2}E_{n,q}\right)  \\
%&=& \int_{\Bbb Z_p} (x+2)^n d\mu_{-1}(x)  \\
&=& 2+ \frac1{q}E_{n,q}=2+\frac1{q} \int_{\Bbb Z_p} x^n q^xd\mu_{-1}(x), \textrm{ for } n>0.\notag
\end{eqnarray*}
Therefore, we obtain the following theorem.
\begin{theorem}
For $n \in \Bbb N$, we have
\begin{eqnarray*}
\int_{\Bbb Z_p} q^{-x}(1-x)^n d\mu_{-1}(x)
= 2+ \frac1{q}\int_{\Bbb Z_p} x^nq^x d\mu_{-1}(x).
\end{eqnarray*}
\end{theorem}
By using Theorem 3, we will study for the $p$-adic integral representation on $\Bbb Z_p$ of the Bernstein polynomials associated with $q$-Euler polynomials in the next section.
\section{ Bernstein polynomials associated with $q$-Euler numbers and polynomials}
Now, we take the $p$-adic integral on $\Bbb Z_p$ for the Bernstein polynomials in (4) as follows:
\begin{eqnarray}
\int_{\Bbb Z_p} B_{k,n}(x)q^x d\mu_{-1}(x)
&=& \int_{\Bbb Z_p} \binom{n}{k} x^k (1-x)^{n-k}q^x d\mu_{-1} (x)\\
&=& \binom{n}{k} \sum_{j=0}^{n-k} \binom{n-k}{j} (-1)^{n-k-j} \int_{\Bbb Z_p} x^{n-j}q^x d\mu_{-1}(x)  \notag\\
&=& \binom{n}{k} \sum_{j=0}^{n-k} \binom{n-k}{j} (-1)^{n-k-j} E_{n-j,q} \notag \\
&=& \binom{n}{k} \sum_{j=0}^{n-k} \binom{n-k}{j} (-1)^{j} E_{k+j,q},~\textrm{ where } n,k\in\Bbb Z_{+}.  \notag
\end{eqnarray}
By the definition of Bernstein polynomials, we see that 
\begin{eqnarray}
B_{k,n}(x)=B_{n-k,n}(1-x),~\textrm{ where } n,k\in\Bbb Z_{+}. 
\end{eqnarray}
Let $n,k\in\Bbb Z_{+}$ with $n>k$. Then, by (15), we get 
\begin{eqnarray*} 
\int_{\Bbb Z_p} q^xB_{k,n}(x) d\mu_{-1}(x) &=& \int_{\Bbb Z_p}q^x B_{n-k,n}(1-x) d\mu_{-1}(x) \notag \\ 
&=& \binom{n}{n-k} \displaystyle\sum_{j=0}^{k} \binom{k}{j} (-1)^{k-j} \int_{\Bbb Z_p} (1-x)^{n-j} q^xd\mu_{-1}(x) \notag \\ 
&=& \binom{n}{k}\displaystyle\sum_{j=0}^{k} \binom{k}{j} (-1)^{k-j} \left(2+ q\int_{\Bbb Z_p} x^{n-j}q^xd\mu_{-1}(x)\right)\notag \\
&=& \binom{n}{k} \displaystyle\sum_{j=0}^{k} \binom{k}{j} (-1)^{k-j} \left(2+qE_{n-j,q^{-1}}\right)  \\
 &=& \left\{
\begin{array}{ll} 2+qE_{n,q^{-1}} \ \ &\hbox{if}\ \ k=0,
\vspace{2mm}\\
\binom{n}{k}q  \displaystyle\sum_{j=0}^{k} \binom{k}{j} (-1)^{k-j}  E_{n-j,q^{-1}} \ \
&\hbox{if}\ \ k>0.
\end{array}\right.  \notag
\end{eqnarray*}
Thus, we obtain the following theorem.

\begin{theorem}
For $n, k \in \Bbb Z_+$ with $n >k$, we have
\begin{eqnarray*}
\int_{\Bbb Z_p}q^{1-x} B_{k, n}(x) d\mu_{-1}(x) =\left\{
\begin{array}{ll} 2q+E_{n,q} \ \ &\hbox{if}\ \ k=0,
\vspace{2mm}\\
\binom{n}{k} \displaystyle\sum_{j=0}^{k} \binom{k}{j} (-1)^{k-j}  E_{n-j,q} \ \
&\hbox{if}\ \ k>0.
\end{array}\right.
\end{eqnarray*}
\end{theorem}

\medskip

By (14) and Theorem 4, we get the following corollary.

\begin{corollary}
For $n, k \in \Bbb Z_+$ with $n >k$, we have
\begin{eqnarray*}
\displaystyle\sum_{j=0}^{n-k} \binom{n-k}{j} (-1)^{j}  E_{k+j,q^{-1}} =\left\{
\begin{array}{ll} 2+\frac1{q}E_{n,q} \ \ &\hbox{if}\ \ k=0,
\vspace{2mm}\\
\displaystyle\sum_{j=0}^{k} \binom{k}{j} (-1)^{k-j} \frac1{q} E_{n-j,q} \ \ &\hbox{if}\ \
k>0.
\end{array}\right.
\end{eqnarray*}
\end{corollary}
%\medskip
For $m, n, k \in \Bbb Z_+$ with $m+n>2k$. Then we get 
\begin{eqnarray*}
& &\int_{\Bbb Z_p} B_{k, n}(x)B_{k, m}(x)q^{-x} d\mu_{-1}(x)
\\  & & \qquad =
\binom{n}{k}\binom{m}{k} \displaystyle\sum_{j=0}^{2k}
\binom{2k}{j}(-1)^{j+2k}\int_{\Bbb Z_p} q^{-x}(1-x)^{n+m-j}d \mu_{-1} (x)\\ 
& & \qquad =\binom{n}{k}\binom{m}{k} \displaystyle\sum_{j=0}^{2k}
\binom{2k}{j}(-1)^{j+2k}q\int_{\Bbb Z_p} (x+2)^{n+m-j}q^xd \mu_{-1} (x)\\ 
& & \qquad =\binom{n}{k}\binom{m}{k} \displaystyle\sum_{j=0}^{2k}
\binom{2k}{j}(-1)^{j+2k}q\left(\frac2{q}+\frac1{q^2} E_{n+m-j,q}\right) \\
& & \qquad =\left\{
\begin{array}{ll} 2+\frac1{q}E_{n+m,q} \ \ &\hbox{if}\ \ k=0,
\vspace{2mm}\\
\binom{n}{k}\binom{m}{k}\displaystyle\sum_{j=0}^{2k} \binom{2k}{j} (-1)^{j+2k}
\frac1{q}E_{n+m-j,q} \ \ &\hbox{if}\ \ k>0.
\end{array}\right.
\end{eqnarray*}
Therefore, we obtain the following theorem.
\begin{theorem}
For $m, n, k \in \Bbb Z_+$ with $m+n>2k$, we have
\begin{eqnarray*}
& &\int_{\Bbb Z_p} B_{k, n}(x)B_{k, m}(x) q^{1-x}d\mu_{-1}(x)
\\& & \qquad =\left\{
\begin{array}{ll} 2q+E_{n+m,q}\ \ &\hbox{if}\ \ k=0,
\vspace{2mm}\\
\binom{n}{k}\binom{m}{k}\displaystyle\sum_{j=0}^{2k} \binom{2k}{j} (-1)^{j+2k}
E_{n+m-j,q} \ \ &\hbox{if}\ \ k>0.
\end{array}\right.
\end{eqnarray*}
\end{theorem}
%\medskip
By using binomial theorem, for $m, n, k \in \Bbb Z_+$, we get
\begin{eqnarray}
& &\int_{\Bbb Z_p} B_{k, n}(x)B_{k, m}(x)q^{1-x}d d\mu_{-1}(x)
\\ & & \qquad =\binom{n}{k}\binom{m}{k} \displaystyle\sum_{j=0}^{n+m-2k}
\binom{n+m-2k}{j}(-1)^{j}\int_{\Bbb Z_p}
x^{j+2k}q^{1-x}d \mu_{-1} (x)\\
& & \qquad =q\binom{n}{k}\binom{m}{k} \displaystyle\sum_{j=0}^{n+m-2k}
\binom{n+m-2k}{j}(-1)^{j} E_{j+2k,q^{-1}}.
\end{eqnarray}
By comparing the coefficients on the both sides of (16) and Theorem 6, we obtain the following corollary. 
\begin{corollary}
Let $m, n, k \in \Bbb Z_+$ with $m+n>2k$. Then we get
\begin{eqnarray*}
& &\displaystyle\sum_{j=0}^{n+m-2k} \binom{n+m-2k}{j}(-1)^{j} E_{j+2k,q^{-1}}
\\& & \qquad =\left\{
\begin{array}{ll} 2+\frac1{q}E_{n+m,q}\ \ &\hbox{if}\ \ k=0,
\vspace{2mm}\\
\frac1{q}\displaystyle\sum_{j=0}^{2k} \binom{2k}{j} (-1)^{j+2k} E_{n+m-j,q} \ \ &\hbox{if}\
\ k>0.
\end{array}\right.
\end{eqnarray*}
\end{corollary}
\medskip
For $s \in \Bbb N$, let $n_1, n_2, \ldots, n_s, k \in \Bbb Z_+$ with $n_1+n_2+\cdots+ n_s>sk.$ By induction, we get 
\begin{eqnarray*}
\int_{\Bbb Z_p} %\prod_{i=1}^{s} 
B_{k, n_1}(x)\cdots   B_{k, n_s}(x)  q^{-x}d\mu_{-1}(x)=
\left(\displaystyle\prod_{i=1}^{s} \binom{n_i}{k} \right)\int_{\Bbb
  Z_p} x^{sk}(1-x)^{n_1+...+n_s}q^{-x}d\mu_{-1}(x)
\end{eqnarray*}
\begin{eqnarray*}
&=&\left(\prod_{i=1}^{s} \binom{n_i}{k} \right) \displaystyle\sum_{j=0}^{sk}\binom{sk}{j} (-1)^{sk+j} \int_{\Bbb Z_p} (1-x)^{n_1+\cdots+n_s -j} q^{-x}d \mu_{-1}(x)
\\ 
&=&\left(\prod_{i=1}^{s} \binom{n_i}{k} \right) \displaystyle\sum_{j=0}^{sk}\binom{sk}{j} (-1)^{sk+j} q\int_{\Bbb Z_p} (x+2)^{n_1+\cdots+n_s -j} q^{x}d \mu_{-1}(x)
\\  &=& \left(\prod_{i=1}^{s} \binom{n_i}{k} \right) \displaystyle\sum_{j=0}^{sk}
\binom{sk}{j}(-1)^{sk+j}q\left(\frac2{q}+\frac1{q^2}E_{n_1+\cdots+n_s-j,q}\right)
\\  &=&\left\{
\begin{array}{ll}  2+\frac1{q}E_{n_1+ \cdots+ n_s,q}\ \ &\hbox{if}\ \ k=0,
\vspace{2mm}\\
\left(\prod_{i=1}^{s} \binom{n_i}{k} \right) \frac1{q}\displaystyle\sum_{j=0}^{sk} \binom{sk}{j}(-1)^{sk+j}E_{n_1+\cdots+n_s -j,q} \ \ &\hbox{if}\ \ k>0.
\end{array}\right.
\end{eqnarray*}
Therefore we obtain the following theorem.
\begin{theorem} Let $s \in \Bbb N$. For $n_1, n_2, \ldots, n_s, k \in
  \Bbb Z_+$ with $n_1+n_2+\cdots+ n_s>sk$, we have
\begin{eqnarray*} & &\int_{\Bbb Z_p} \left(\prod_{i=1}^{s} B_{k, n_i}(x) \right) q^{1-x}d\mu_{-1}(x) \\& & \qquad =\left\{ \begin{array}{ll} 2q+E_{n_1+n_2+\cdots+ n_s,q}\ \ &\hbox{ if }\ \ k=0, \vspace{2mm}\\ \left(\displaystyle\prod_{i=1}^{s} \binom{n_i}{k} \right)\displaystyle\displaystyle\sum_{j=0}^{sk} \binom{s k}{j}(-1)^{sk+j} E_{n_1+n_2+\cdots+ n_s-j,q} \ \ &\hbox{ if }  \ \ k>0.
\end{array}\right. \end{eqnarray*} 
\end{theorem}

For $n_1, n_2, \ldots, n_s, k \in \Bbb Z_+$ by binomial theorem, we get 
\begin{eqnarray} 
& &
\int_{\Bbb Z_p} \left(\prod_{i=1}^{s} B_{k, n_i}(x) \right) q^{-x}d\mu_{-1}(x) \\
& & = \binom{n_1}{k}...\binom{n_s}{k}\displaystyle\sum_{j=0}^{n_1 + \cdots +n_s - sk} \binom{n_1 + \cdots +n_s -s k}{j}(-1)^j \int_{\Bbb Z_p} x^{j+sk}q^{-x}d\mu_{-1}(x) \notag\\
& & = \binom{n_1}{k}...\binom{n_s}{k}\displaystyle\sum_{j=0}^{n_1 + \cdots +n_s - sk} \binom{n_1 + \cdots +n_s -s k}{j}(-1)^jE_{j+s k,q^{-1}}.\notag
 \end{eqnarray} 

By using (19) and Theorem 8, we obtain the following corollary.

\begin{corollary} Let $s \in \Bbb N$.
For $n_1, n_2, \ldots, n_s, k \in \Bbb Z_+$ with $n_1+n_2+\cdots+
n_s>sk$, we have
\begin{eqnarray*}
& &  \displaystyle\sum_{j=0}^{n_1 + \cdots +n_s - sk} \binom{n_1 + \cdots +n_s -
s k}{j}(-1)^{j} E_{j+s k,q^{-1}}
\\& & \qquad =\left\{
\begin{array}{ll}  2+\frac1{q}E_{n_1+n_2+\cdots+
n_s,q}\ \ &\hbox{if}\ \ k=0,
\vspace{2mm}\\
\frac1{q}\displaystyle\sum_{j=0}^{sk} \binom{s k}{j}(-1)^{sk+j} E_{n_1+n_2+\cdots+ n_s-j,q}
\ \ &\hbox{if}\ \ k>0.
\end{array}\right.
\end{eqnarray*}
\end{corollary}

\end{document}